\documentclass[aop,MSNbibl,citesort,dvips]{arximspdf}
\usepackage{mathbh}

\doi{10.1214/10-AOP540}
\volume{38}
\issue{6}
\pubyear{2010}
\firstpage{2295}
\lastpage{2321}

\makeatletter

\newtheorem{theorem}{Theorem}[section]
\newtheorem{proposition}[theorem]{Proposition}
\newtheorem{lemma}[theorem]{Lemma}

\newproclaim{remark}{Remark}

\makeatother

\begin{document}
\begin{frontmatter}

\title{Construction of an Edwards' probability\\ measure on $\mathcal{C}
(\mathbb{R}_+, \mathbb{R})$}
\runtitle{Construction of an Edwards' probability measure}

\begin{aug}
\author[A]{\fnms{Joseph} \snm{Najnudel}\corref{}\ead[label=e1]{joseph.najnudel@math.uzh.ch}}
\runauthor{J. Najnudel}
\affiliation{Universit\"{a}t Z\"{u}rich}
\address[A]{Institut f\"{u}r Mathematik\\
Universit\"{a}t Z\"{u}rich\\
Winterthurer strasse 190\\
CH-8057 Z\"{u}rich\\
Switzerland\\
\printead{e1}} 
\end{aug}

\received{\smonth{1} \syear{2008}}
\revised{\smonth{2} \syear{2010}}

%
\begin{abstract}
In this article, we prove that the measures $\mathbb{Q}_T$ associated
to the one-dimensional Edwards' model on the interval $[0,T]$ converge
to a limit measure $\mathbb{Q}$ when $T$ goes to infinity, in the
following sense: for all $s \geq0$ and for all events $\Lambda_s$
depending on the canonical process only up to time $s$,
$\mathbb{Q}_T (\Lambda_s) \rightarrow\mathbb{Q} (\Lambda_s)$.

Moreover, we prove that, if $\mathbb{P}$ is Wiener measure, there
exists a martingale $(D_s)_{s \in\mathbb{R}_+}$ such that $\mathbb{Q}
(\Lambda_s) = \mathbb{E}_{\mathbb{P}} (\mathbh{1}_{\Lambda_s} D_s)$,
and we give an explicit expression for this martingale.
\end{abstract}

%

\end{frontmatter}

\section{Introduction and statement of the main theorems}

Edwards' model is a model for polymers chains, which is defined by
considering Brownian motion ``penalized'' by the ``quantity'' of its
self-intersections (see also \cite{e65}). More precisely, for $d \in
\mathbb{N}^*$, and $T > 0$, let $\mathbb{P}_T^{(d)}$ be Wiener measure
on the space $\mathcal{C}([0,T], \mathbb{R}^d)$, and $(X_t^{(d)})_{t
\in[0,T]}$ the corresponding canonical process. The $d$-dimensional
Edwards' model on $[0,T]$ is defined by the probability measure
$\mathbb{Q}_T^{(d), \beta}$ on $\mathcal{C}([0,T], \mathbb{R}^d)$
such that, very informally,
%
%
\begin{equation} \label{1}
\mathbb{Q}_T^{(d), \beta} = \frac{\exp( - \beta\int_0^T \int_0^T
\delta(X_s^{(d)} - X_u^{(d)}) \,ds \,du )}
{ \mathbb{P}_T^{(d)} [ \exp( - \beta\int_0^T \int_0^T
\delta(X_s^{(d)} - X_u^{(d)} ) \,ds \,du ) ] }
\cdot\mathbb{P}_T^{(d)},
\end{equation}
where $\beta$ is a strictly positive parameter, and $\delta$ is Dirac
measure at zero.

(In this article, we always denote by $\mathbb{Q} [V]$ the expectation
of a random variable $V$ under the probability $\mathbb{Q}$.)

Of course, (\ref{1}) is not really the definition of a probability
measure, since the integral with respect to Dirac measure is not well
defined. However, it has been proven that one can define rigorously the
measure $\mathbb{Q}_T^{(d), \beta}$ for $d=1, 2, 3$, by giving a
meaning to (\ref{1}) (for $d \geq4$, the Brownian path has no
self-intersection, so the measure $\mathbb{Q}_T^{(d), \beta}$ has to be
equal to $\mathbb{P}_T^{(d)}$).

In particular, for $ d= 1$, one has formally the equality
%
%
\begin{equation} \label{2}
\int_0^T \int_0^T \delta\bigl(X_s^{(1)} - X_u^{(1)}\bigr) \,ds \,du = \int
_{- \infty}^{\infty} (L_T^y)^2 \,dy,
\end{equation}
where $(L_T^y)_{y \in\mathbb{R}}$ is the continuous family of local
times of $(X_s^{(1)})_{s \leq T}$ (which is $ \mathbb
{P}_T^{(1)}$-almost surely well defined).

Therefore, one can take the following (rigorous) definition:
\[
\mathbb{Q}_T^{(1), \beta} = \frac{\exp( - \beta\int_{-\infty
}^{\infty} (L_T^y)^2 \,dy )}
{ \mathbb{P}_T^{(1)} [ \exp( - \beta\int_{-\infty}^{\infty}
(L_T^y)^2 \,dy ) ] }
\cdot\mathbb{P}_T^{(1)}.
\]
Under $\mathbb{Q}_T^{(1), \beta}$, the canonical process has a
ballistic behavior; more precisely, Westwater (see \cite{w85}) has
proven that for $T \rightarrow\infty$, the law of $\frac{
X_T^{(1)}}{T}$ under $\mathbb{Q}_T^{(1), \beta}$ tends to $\frac{1}{2}
(\delta_{b^* \beta^{1/3}}
+ \delta_{-b^* \beta^{1/3}} )$, where $\delta_x$ is Dirac measure at
$x$ and $b^*$ is a universal constant (approximately equal to 1.1).

This result was improved in \cite{hhk97} (see also \cite{h98}), where
van der Hofstad, den Hollander and K{\"o}nig show that $\frac{|
X_T^{(1)} | - b^* \beta^{1/3} T}{\sqrt{T}}$ tends in law to a centered
Gaussian variable, which has a variance equal to a universal constant
(approximately equal to 0.4; in particular, smaller than one).

Moreover, in \cite{hhk03}, the authors prove large deviation results
for the variable $X_T$ under $\mathbb{Q}_T^{(1), \beta}$.

In dimension $2$, the problem of the definition of Edwards' model was
solved by Varadhan (see \cite{v69,lg85,r86}). In this
case, it is possible to give a rigorous definition of $I := \int_0^T
\int_0^T \delta(X_s^{(2)} - X_u^{(2)}) \,ds \,du$, but this
quantity appears to be equal to infinity. However, if one formally
subtracts its expectation (i.e., one considers the quantity:
$I - \mathbb{P}_T^{(2)} [I]$), one can define a finite random variable
which has negative exponential moments of any order; therefore, if we
replace $\int_0^T \int_0^T \delta(X_s^{(2)} - X_u^{(2)}) \,ds
\,du$ by this random variable in (\ref{1}), we obtain a rigorous
definition of $\mathbb{Q}_T^{(2), \beta}$. Moreover, this probability
is absolutely continuous with respect to Wiener measure.

In dimension 3 (the most difficult case), subtracting the expectation
(this technique is also called ``Varadhan renormalization'') is not
sufficient to define Edwards' model. However, by a long and difficult
construction, Weswater (see \cite{w80,w82}) has proven that it
is possible to define the probability $\mathbb{Q}_T^{(3), \beta}$;
this construction has been simplified by Bolthausen in \cite{b93} (at
least if $\beta$ is small enough).
Moreover, the measures $(\mathbb{Q}_T^{(3), \beta})_{\beta\in
\mathbb
{R}_+^*}$ are mutually singular, and singular with respect to Wiener
measure.

The behavior of the canonical process under $\mathbb{Q}_T^{(d), \beta
}$, as $T \rightarrow\infty$, is essentially unknown for $d = 2$ and
$d=3$. One conjectures that the following convergence holds:
\[
\mathbb{Q}_T^{(3), \beta} [\Vert X_T\Vert] \rightarrow D T^{\nu},
\]
where $D > 0$ depends only on $d$ and $\beta$, and where $\nu$ is equal
to $3/4$ for $d=2$ and approximately equal to $0.588$ for $d=3$ (see
\cite{h98}, Chapter 1).

At this point, we note that all the measures considered above are
defined on finite interval trajectories [exactly, on $\mathcal{C}
([0,T], \mathbb{R})$].

An interesting question is the following: is it possible to define
Edwards' model on trajectories indexed by $\mathbb{R}_+$?

More precisely, if $\mathbb{P}^{(d)}$ is Wiener measure on $\mathcal
{C}(\mathbb{R}_+, \mathbb{R}^d)$ and $(X_s^{(d)})_{s \in\mathbb{R}_+}$
the corresponding canonical process, is it possible to define a measure
$\mathbb{Q}^{(d), \beta}$ (for all $\beta> 0$) such that, informally,
\[
\mathbb{Q}^{(d), \beta} = \frac{\exp( - \beta\int_0^{\infty}
\int_0^{\infty} \delta(X_s^{(d)} -
X_u^{(d)}) \,ds \,du )}
{ \mathbb{P}^{(d)} [ \exp( - \beta\int_0^{\infty} \int
_0^{\infty} \delta(X_s^{(d)} - X_u^{(d)} ) \,ds \,du )
] }
\cdot\mathbb{P}^{(d)} ?
\]
In this article, we give a positive answer to this question in
dimension one. The construction of the corresponding measure is
analogous to the construction given by Roynette, Vallois and Yor in
their articles about penalisation (see \cite{rvy06j,rvy05,rvy061,rvy062}).

More precisely, let us replace the
notation $\mathbb{P}^{(1)}$ by $\mathbb{P}$ for the standard Wiener
measure and the notation $(X_s^{(1)})_{s \in\mathbb{R}_+}$ by
$(X_s)_{s \in\mathbb{R}_+}$ for the canonical process. If $(\mathcal
{F}_s)_{s \in\mathbb{R}_+}$ is the natural filtration of $X$, and if
for all $T \in\mathbb{R}_+$, the measure $\mathbb{Q}_T^{\beta}$ is
defined by
\[
\mathbb{Q}_T^{\beta} = \frac{\exp( - \beta\int_{-\infty}^{\infty
} (L_T^y)^2 \,dy )}
{ \mathbb{P} [ \exp( - \beta\int_{-\infty}^{\infty}
(L_T^y)^2 \,dy ) ] }
\cdot\mathbb{P},
\]
where $(L_T^y)_{T \in\mathbb{R}_+, y \in\mathbb{R}}$ is the jointly
continuous version of the local times of $X$ ($\mathbb{P}$-almost
surely well defined), the following theorem holds.
\begin{theorem} \label{main}
For all $\beta> 0$, there exists a unique probability measure $\mathbb
{Q}^{\beta}$ such that for all $s \geq0$, and for all events $\Lambda
_s \in\mathcal{F}_s$,
%
%
\begin{equation}
\mathbb{Q}_T^{\beta} (\Lambda_s) \mathop{\longrightarrow}\limits
_{ T \rightarrow\infty}
\mathbb{Q}^{\beta} (\Lambda_s).
\end{equation}
\end{theorem}

Theorem \ref{main} is the main result of our article.

Let us remark that if $\Lambda_s \in\mathcal{F}_s$ ($s \geq0$) and
$\mathbb{P} (\Lambda_s) = 0$, then $\mathbb{Q}_T^{\beta} (\Lambda
_s) =
0$, since $\mathbb{Q}_T^{\beta}$ is, by definition, absolutely
continuous with respect to $\mathbb{P}$. Hence, if Theorem~\ref{main}
is assumed, $\mathbb{Q}^{\beta} (\Lambda_s)$ is equal to zero.

Therefore, the restriction of $\mathbb{Q}^{\beta}$ to $\mathcal{F}_s$
is absolutely continuous with respect to the restriction of $\mathbb
{P}$ to $\mathcal{F}_s$, and there exists a $\mathbb{P}$-martingale
$(D^{\beta}_s)_{s \geq0}$ such that, for all~$s$,
\[
\mathbb{Q}^{\beta} _{| \mathcal{F}_s} = D^{\beta}_s\cdot\mathbb{P}
_{| \mathcal{F}_s}.
\]
In our proof of Theorem \ref{main}, we obtain an explicit formula for
the martingale $(D^{\beta}_s)_{s \geq0}$. However, we need to define
other notation before giving this formula.

Let $\nu$ be the measure on $\mathbb{R}_+^*$, defined by $\nu(dx)= x
\,dx$, and let $L^{2} (\nu)$ be the set of functions $g$ from $\mathbb
{R}_+^*$ to $\mathbb{R}$ such that
\[
\int_{0}^{\infty} [g(x)]^2 \nu(dx) < \infty,
\]
equipped with the scalar product
\[
\langle g|h \rangle= \int_0^{\infty} g(x) h(x) \nu(dx).
\]
The operator $\mathcal{K}$ defined from $L^{2}(\nu) \cap\mathcal{C}^2
(\mathbb{R}_+^*)$ to $\mathcal{C} (\mathbb{R}_+^*)$ by
%
%
\begin{equation} \label{propre}
[\mathcal{K} (g)] (x) = 2 g''(x) + \frac{2 g'(x)}{x} - x g(x)
\end{equation}
is the infinitesimal generator of the process $2R$ killed at rate $x$
at level $x$, where $R$ is a Bessel process of dimension two; it is a
Sturm--Liouville operator, and there exists an orthonormal basis
$(e_n)_{n \in\mathbb{N}}$ of $ L^2(\nu) $, consisting of
eigenfunctions of $\mathcal{K}$, with the corresponding negative
eigenvalues: $- \rho_0 > - \rho_1 \geq- \rho_2 \geq- \rho_3 \geq
\cdots,$ where $\rho:= \rho_0$ is in the interval $[2.18,2.19]$.

Moreover, the functions $(e_n)_{n \in\mathbb{N}}$ are analytic and
bounded (they tend to zero at infinity, faster than exponentially), and
$e_0$ is strictly positive (these properties are quite classical, and
they are essentially proven in \cite{h98}, Chapters 2 and 3; see also~\cite{hlac}).

Now, for $l \in\mathbb{R}_+$, let us denote by $(Y_l^y)_{y \in
\mathbb
{R}}$ a process from $\mathbb{R}$ to $\mathbb{R}_+$ such that:
\begin{itemize}
\item$(Y_l^{-y})_{y \geq0}$ is a squared Bessel process of dimension
zero, starting at $l$.
\item$(Y_l^{y})_{y \geq0}$ is an independent squared Bessel process
of dimension two.
\end{itemize}
Moreover, let $f$ be a continuous function with compact support from
$\mathbb{R}$ to $\mathbb{R}_+$, and let $M$ be a strictly positive real
such that $f(x) = 0$ for all $x \notin[-M,M]$. We define the following
quantities:
\begin{eqnarray*}
A_+^{\beta, M} (f) &=& \int_0^{\infty} dl\, \mathbb{E} \bigl[ e^{ \int
_{-\infty}^M [ - \beta(Y_l^y + f(y))^2
+ \rho\beta^{2/3} Y_l^y ] \,dy } e_0 (\beta^{1/3} Y_l^M )
\bigr],
\\
A_-^{\beta, M} (f) &=& A_+^{\beta, M} (\widetilde{f}),
\end{eqnarray*}
where $\widetilde{f}$ is defined by $\widetilde{f} (x) = f(-x)$, and
\[
A^{\beta, M} (f) = A_+^{\beta, M} (f)+ A_-^{\beta, M} (f).
\]
With this notation, we can state the following theorem, which gives
an explicit formula for the martingale $(D^{\beta}_s)_{s \geq0}$.
\begin{theorem} \label{mart}
For all $\beta>0$ and for all continuous and positive functions $f$
with compact support, the quantity $A^{\beta,M} (f)$ is finite,
different from zero, and does not depend on the choice of $M > 0$ such
that $f=0$ outside the interval $[-M,M]$; therefore, we can write:
$A^{\beta} (f) := A^{\beta,M} (f)$. Moreover, for all $s \geq0$, the
density $D^{\beta}_s$ of the restriction of $\mathbb{Q}^{\beta}$ to
$\mathcal{F}_s$, with respect to the restriction of $\mathbb{P}$ to
$\mathcal{F}_s$, is given by the equality
%
%
\begin{equation} \label{00}
D^{\beta}_s = e^{\rho\beta^{2/3} s}\cdot\frac{A^{\beta}
(L_s^{\bullet+ X_s})} { A^{\beta} (0)},
\end{equation}
where $L_s^{\bullet+ X_s}$ denotes the function $F$ [which depends on
the trajectory $(X_u)_{u \leq s}$] such that $F(y) = L_s^{y + X_s}$ for
all $y \in\mathbb{R}$.
\end{theorem}
\begin{remark*} The independence of $A^{\beta,M} (f)$ with respect to
$M$ (provided the support of $f$ is included in $[-M,M]$) can be
checked directly by using the fact that
%
%
\begin{equation} \biggl( \exp\biggl( \int_0^x [ - \beta(Y_l^y)^2 +
\rho\beta^{2/3} Y_l^y ] \,dy \biggr) e_0 (\beta^{1/3} Y_l^x
) \biggr)_{x \geq0}
\end{equation}
is a martingale, property which can be easily proven by using the
differential equation
satisfied by $e_0$.
\end{remark*}

For $l > 0$, $\mu\in\mathbb{R}$ and $v > 0$, let us now define the
following quantity:
%
%
\begin{equation} \label{kim28}
K_l^{(\mu)} (v) = \alpha_l(v) e^{\mu v} \mathbb{E} [ e^{-2 \int
_0^v V_u^{(l/2,v)} \,du } ],
\end{equation}
where $\alpha_l (v) = \frac{l}{\sqrt{8 \pi v^3}} e^{-l^2/8v}$
denotes the density of the first hitting time of zero of a Brownian
motion starting at $l/2$ (or equivalently, the density of the last
hitting time of $l/2$ of a standard Bessel process of dimension 3),
and $ (V_u^{(l/2,v)})_{u \leq v}$ is the bridge of a Bessel process of
dimension 3 on $[0,v]$, starting at $l/2$ and ending at $0$.

To simplify the notation, we set
\[
K_l (v) = K_l^{(0)} (v).
\]
Moreover, let us consider, for $v > 0$, the function $\chi_v$ defined by
%
%
\begin{equation} \label{guet}
\chi_v(l) =\frac{K_l (v)}{l} = \frac{1}{\sqrt{8 \pi
v^3}}e^{-l^2/8v} \mathbb{E} [ e^{-2 \int_0^v V_u^{(l/2,v)}
\,du } ]
\end{equation}
for $l > 0$.

With this notation, Theorem \ref{mart} is a essentially a consequence
of the two
propositions stated below.
\begin{proposition} \label{prop1}
When $T$ goes to infinity,
%
%
\begin{equation}\label{b}
e^{\rho T} \mathbb{P} \bigl[ e^{- \int_{- \infty}^{\infty} [L_T^y +
f(y)]^2 \,dy} \mathbh{1}_{X_T \in[0,M]} \bigr] \longrightarrow0.
\end{equation}
\end{proposition}
\begin{proposition} \label{prop2}
When $T$ goes to infinity,
%
%
\begin{equation} \label{a}
e^{\rho T} \mathbb{P} \bigl[ e^{- \int_{- \infty}^{\infty} [L_T^y +
f(y)]^2 \,dy} \mathbh{1}_{X_T \geq M} \bigr] \longrightarrow K
A_+^{1,M} (f) < \infty,
\end{equation}
where $K \in\mathbb{R}_+^*$ is a universal constant (in particular,
$K$ does not depend on $f$ and~$M$).

Moreover, for all $ v> 0$, $\chi_v \in L^2 (\nu)$ and the constant $K$
is given by
the formula
\[
K = \int_0^{\infty} e^{\rho v} \langle\chi_v | e_0 \rangle \,dv
< \infty.
\]
\end{proposition}

In the proof of these two propositions, we use essentially the same
tools as in the papers
by van der Hofstad, den Hollander and K\"onig. In particular, for $f =
0$, Propositions \ref{prop1} and \ref{prop2} are consequences of
Proposition 1 of \cite{hhk97}.

However, for a general function $f$, it is not obvious that one can
deduce directly our results from the material of
\cite{hhk97} and \cite{hhk03}, since for $X_T > 0$, one has to deal
with the family of local times of the canonical process on the intervals
$\mathbb{R}_-$, $[0, X_T]$ and $[X_T, \infty)$ as for $f=0$, but also
on the support of $f$. Moreover, some typos in \cite{hhk97} make the
argument as written incorrect. For this reason, we present a proof of
this result in a different way than was done in \cite{hhk97}.

The next sections of this article are organized as follows. In Section
\ref{sec2}, we prove that Propositions \ref{prop1} and \ref{prop2} imply
Theorems \ref{main} and \ref{mart}; in Section \ref{sec3}, we prove Proposition
\ref{prop1}. The proof of Proposition \ref{prop2} is split into two
parts: the first
one is given in Section \ref{sec4}; the second one, for which one needs some
estimates of different quantities, is given in Section \ref{sec6}, after the
proof of these estimates in Section \ref{sec5}. In Section \ref{sec7}, we make a
conjecture on the behavior of the canonical process under the limit
measure $\mathbb{Q}^{\beta}$.

\section{\texorpdfstring{Proof of Theorems \protect\ref{main} and \protect\ref{mart} by assuming
Propositions \protect\ref{prop1} and \protect\ref{prop2}}{Proof of Theorems 1.1 and 1.2 by assuming
Propositions 1.3 and 1.4}}\label{sec2}

Let us begin to prove the following result, which is essentially a
consequence of Brownian scaling.
\begin{proposition} \label{18}
Let us assume Propositions \ref{prop1} and \ref{prop2}. For any
positive continuous function $f$ with compact support included in
$[-M,M]$, and for all $\beta> 0$,
%
%
\begin{equation} \label{123}
e^{\rho\beta^{2/3} T} \mathbb{P} \bigl[ e^{- \beta\int_{-\infty
}^{\infty} [L_T^y + f(y)] ^2 \,dy } \bigr] \longrightarrow K \beta
^{1/3} A^{\beta, M} (f) < \infty,
\end{equation}
when $T$ goes to infinity.
\end{proposition}
\begin{pf}
Propositions \ref{prop1} and \ref{prop2} imply
\[
e^{\rho T} \mathbb{P} \bigl[ e^{- \int_{-\infty}^{\infty} [L_T^y +
f(y)]^2 \,dy } \mathbh{1}_{X_T \geq0} \bigr]
\mathop{\longrightarrow}\limits_{T \rightarrow\infty}
K A_{+}^{1,M} (f) < \infty.
\]
Now, $((L_T^{-y})_{y \in\mathbb{R}}, - X_T)$ and $((L_T^{y})_{y \in
\mathbb{R}}, X_T)$ have the same law; hence,
\begin{eqnarray*}
&&
e^{\rho T} \mathbb{P} \bigl[ e^{- \int_{-\infty}^{\infty} [L_T^y +
f(y)]^2 \,dy } \mathbh{1}_{X_T \leq0} \bigr] \\
&&\qquad = e^{\rho T}
\mathbb{P} \bigl[ e^{- \int_{-\infty}^{\infty} [L_T^{-y} + f(y)]^2
\,dy } \mathbh{1}_{-X_T \leq0} \bigr] \\
&&\qquad = e^{\rho T}
\mathbb{P} \bigl[ e^{- \int_{-\infty}^{\infty} [L_T^y + f(-y)]^2 \,dy }
\mathbh{1}_{X_T \geq0} \bigr]
\mathop{\longrightarrow}\limits_{T \rightarrow\infty}
K A_{+}^{1,M} (\widetilde{f})\\
&&\qquad =
K A_{-}^{1,M} (f),
\end{eqnarray*}
which is finite.

Therefore,
\[
e^{\rho T} \mathbb{P} \bigl[ e^{- \int_{-\infty}^{\infty} [L_T^y +
f(y)]^2 \,dy } \bigr]
\mathop{\longrightarrow}\limits_{T \rightarrow\infty}
K A^{1,M} (f) < \infty.
\]
Now, let us set: $\alpha= \beta^{1/3}$.
By Brownian scaling, $(L_{T \alpha^2}^{y \alpha})_{y \in\mathbb{R}} $
and $(\alpha L_T^y)_{y \in\mathbb{R}}$ have the same law. Consequently,
\begin{eqnarray*}
&&
e^{\rho\alpha^2 T} \mathbb{P} \bigl[ e^{- \beta\int_{-\infty
}^{\infty} [L_T^y + f(y)]^2 \,dy } \bigr] \\
&&\qquad = e^{\rho\alpha^2 T}
\mathbb{P} \bigl[ e^{- \alpha\int_{-\infty}^{\infty} [L_{T \alpha
^{2}}^{y \alpha} + \alpha f(y)]^2 \,dy } \bigr] \\
&&\qquad =
e^{\rho\alpha^{2} T} \mathbb{P} \bigl[ e^{- \int_{-\infty}^{\infty}
[L_{T \alpha^2}^{z} + \alpha f(z \alpha^{-1})]^2 \,dz } \bigr]
\mathop{\longrightarrow}\limits_{T \rightarrow\infty}
K A^{1,M \alpha} (f_{\alpha}) \\
&&\qquad< \infty,
\end{eqnarray*}
where $f_{\alpha}$, defined by $f_{\alpha}(z) = \alpha f(z \alpha
^{-1})$, has a support included in $[-M \alpha, M \alpha]$.

Therefore, Proposition \ref{18} is proven if we show that $A^{1,M
\alpha
} (f_{\alpha}) = \alpha A^{\beta, M} (f)$.

Now, by change of variable and scaling property of squared Bessel
processes,
%
%
\begin{eqnarray} \label{1111}\quad
A_+^{1,M \alpha} (f_{\alpha}) &=& \int_0^{\infty} dl\, \mathbb{E}
\bigl[e ^ { \int_{-\infty}^{M \alpha}
[ - ( Y_l^y + \alpha f(y \alpha^{-1}) )^2 + \rho Y_l^y ] \,dy}
e_0 (Y_l^{M \alpha}) \bigr]
\nonumber\\
&=& \int_0^{\infty} dl\, \mathbb{E} \bigl[e ^ { \alpha
\int_{-\infty}^{M}
[ - ( Y_l^{z \alpha} + \alpha f(z) )^2 + \rho Y_l^{z \alpha}
] \,dz} e_0 (Y_l^{M \alpha}) \bigr]
\nonumber\\
&=& \int_0^{\infty} dl\, \mathbb{E} \bigl[e ^ { \beta
\int_{-\infty}^{M}
[ - ( Y_{l \alpha^{-1}}^{z} + f(z) )^2 + \rho\alpha^{-1} Y_{l
\alpha^{-1}} ^{z} ] \,dz} e_0 ( \alpha Y_{l \alpha^{-1}}^{M})\bigr]
\\
&=& \alpha\int_0^{\infty} dl\, \mathbb{E} \bigl[e ^ {
\beta\int_{-\infty}^{M}
[ - ( Y_{l}^{z} + f(z) )^2 + \rho\alpha^{-1} Y_{l} ^{z} ]
\,dz} e_0 ( \alpha Y_{l}^{M}) \bigr]\nonumber\\
&=& \alpha A_+^{\beta, M} (f).\nonumber
\end{eqnarray}
By replacing $f$ by $\widetilde{f}$, one obtains
%
%
\begin{equation} \label{1112}
A_-^{1,M \alpha} (f_{\alpha}) = \alpha A_-^{\beta, M} (f),
\end{equation}
and by adding (\ref{1111}) and (\ref{1112}),
\[
A^{1,M \alpha} (f_{\alpha}) = \alpha A^{\beta, M} (f),
\]
which proves Proposition \ref{18}.
\end{pf}

At this point, we remark that $A^{\beta} (f) := A^{\beta,M}(f)$ does
not depend on $M$ (as written in Theorem \ref{mart}), since $M$ does
not appear in the left-hand side of (\ref{123}).

Now, let $T > s$ be in $\mathbb{R}_+$. One has, for all $y
\in\mathbb{R}$,
\[
L_T^y = L_s^y + \widetilde{L}_{T-s}^{y-X_s},
\]
where $\widetilde{L}$ is the continuous family of local times of the
process $(X_{s+u}-X_s)_{u \geq0}$.

Therefore, for all $\beta> 0$,
\begin{eqnarray*}
\mathbb{P} \bigl[ e^{- \beta\int_{-\infty}^{\infty} (L_T^y)^2 \,dy} |
\mathcal{F}_s \bigr]
&=& \mathbb{P} \bigl[ e^{- \beta\int_{-\infty}^{\infty} (L_s^y +
\widetilde{L}_{T-s}^{y-X_s})^2 \,dy} | \mathcal{F}_s \bigr]
\\
&=& \mathbb{P} \bigl[ e^{- \beta\int_{-\infty}^{\infty}
(L_s^{y+X_s} + \widetilde{L}_{T-s}^{y})^2 \,dy} | \mathcal{F}_s
\bigr].
\end{eqnarray*}
Under $\mathbb{P}$ and conditionally on $\mathcal{F}_s$,
$(L_s^{y+X_s})_{y \in\mathbb{R}}$ is fixed and by Markov property, $
(X_{s+u}-X_s)_{u \geq0}$ is a standard Brownian motion.

Hence, if we assume Propositions \ref{prop1} and \ref{prop2}, we
obtain, by using Proposition~\ref{18},
\[
e^{\rho(T-s) \alpha^2} \mathbb{P} \bigl[ e^{- \beta\int_{-\infty
}^{\infty} (L_T^y)^2 \,dy} | \mathcal{F}_s \bigr]
\mathop{\longrightarrow}\limits_{T \rightarrow\infty}
K \alpha A^{\beta} (
L_s^{\bullet+ X_s}).
\]
Moreover,
\begin{eqnarray*}
e^{\rho(T-s) \alpha^2} \mathbb{P} \bigl[ e^{- \beta\int_{-\infty
}^{\infty} (L_T^y)^2 \,dy} | \mathcal{F}_s \bigr] & \leq & e^{\rho
(T-s) \alpha^2} \mathbb{P} \bigl[ e^{- \beta\int_{-\infty}^{\infty} (
\widetilde{L}_{T-s}^{y})^2 \,dy} | \mathcal{F}_s \bigr] \\
& = & e^{\rho(T-s) \alpha^2} \mathbb{P} \bigl[ e^{- \beta\int_{-\infty
}^{\infty} ( L_{T-s}^{y})^2 \,dy} \bigr] \\ & \leq &2 K
\alpha A^{\beta} (0) < \infty,
\end{eqnarray*}
if $T-s$ is large enough.
On the other hand,
\[
e^{ \rho T \alpha^2} \mathbb{P} \bigl[ e^{- \beta\int_{-\infty
}^{\infty} (L_T^y)^2 \,dy} \bigr]
\mathop{\longrightarrow}\limits_{T \rightarrow\infty}
K \alpha A^{\beta} (0),
\]
and for $T$ large enough,
\[
e^{ \rho T \alpha^2} \mathbb{P} \bigl[ e^{- \beta\int_{-\infty
}^{\infty} (L_T^y)^2 \,dy} \bigr] \geq\frac{K}{2} \alpha A^{\beta}
(0).
\]
Now, for all $\beta$ and $f$, $ A^{\beta} (f)$ is different from zero
(as written in Theorem \ref{mart}), since it is the integral of a
strictly positive quantity. Therefore,
\[
\frac{\mathbb{P} [ e^{- \beta\int_{-\infty}^{\infty} (L_T^y)^2
\,dy} | \mathcal{F}_s ]}
{ \mathbb{P} [ e^{- \beta\int_{-\infty}^{\infty} (L_T^y)^2
\,dy} ]}
\mathop{\longrightarrow}\limits_{T \rightarrow\infty}
e^{\rho
\alpha^2 s} \frac{A^{\beta} ( L_s^{\bullet+ X_s})} {A^{\beta}
(0)},
\]
and for $s$ fixed and $T$ large enough
\[
\frac{\mathbb{P} [ e^{- \beta\int_{-\infty}^{\infty} (L_T^y)^2
\,dy} | \mathcal{F}_s ]}
{ \mathbb{P} [ e^{- \beta\int_{-\infty}^{\infty} (L_T^y)^2
\,dy} ]} \leq4 e^{\rho\alpha^2 s}
< \infty.
\]
Consequently, for all $s \geq0$ and $\Lambda_s \in\mathcal{F}_s$, by
dominated convergence
\[
\mathbb{P} \biggl[ \mathbh{1}_{\Lambda_s} \frac{\mathbb{P} [
e^{- \beta\int_{-\infty}^{\infty} (L_T^y)^2 \,dy} | \mathcal{F}_s
]}
{ \mathbb{P} [ e^{- \beta\int_{-\infty}^{\infty} (L_T^y)^2
\,dy} ]} \biggr] \mathop{\longrightarrow}\limits_{T \rightarrow\infty}
\mathbb{P} \biggl[ \mathbh{1}_{\Lambda_s} e^{\rho\alpha^2 s} \frac
{A^{\beta} ( L_s^{\bullet+ X_s})} {A^{\beta} (0)} \biggr].
\]
Hence,
\[
\mathbb{Q}_T^{\beta} (\Lambda_s)
\mathop{\longrightarrow}\limits_{T \rightarrow\infty}
\mathbb{P} (\mathbh{1}_{\Lambda_s} D^{\beta}_s),
\]
where $D^{\beta}_s$ is defined by (\ref{00}):
\[
D_s^{\beta} = e^{\rho\beta^{2/3} s} \cdot \frac{A^{\beta
}(L_s^{\bullet+
X_s})}{A^{\beta}(0)}.
\]
This convergence implies Theorems \ref{main} and \ref{mart}.

\section{\texorpdfstring{Proof of Proposition \protect\ref{prop1}}{Proof of Proposition 1.3}}\label{sec3}

If $f$ is a continuous function from $\mathbb{R}$ to $\mathbb{R}_+$
with compact support included in $[-M,M]$, one has
%
%
\begin{eqnarray} \label{007}\quad
\mathbb{P} \bigl[ e^{- \int_{- \infty}^{\infty} [L_T^y + f(y)]^2 \,dy}
\mathbh{1}_{X_T \in[0,M]} \bigr] &\leq& \mathbb{P} \bigl[ e^{- \int
_{- \infty}^{\infty} (L_T^y)^2 \,dy} \mathbh{1}_{X_T \in[0,M]}
\bigr]
\nonumber\\[-8pt]\\[-8pt]
&=& \mathbb{P} \bigl[ e^{- T^{3/2} \int_{- \infty
}^{\infty} (L_1^y)^2 \,dy} \mathbh{1}_{X_1 \in[0,M/\sqrt{T}]}\bigr]\nonumber
\end{eqnarray}
by scaling properties of Brownian motion.

Hence, the right-hand side of (\ref{007}) is decreasing with $T$, which
implies (for $T > 1$)
\begin{eqnarray*}
&&
e^{\rho T} \mathbb{P} \bigl[ e^{- \int_{- \infty}^{\infty} [L_T^y +
f(y)]^2 \,dy} \mathbh{1}_{X_T \in[0,M]} \bigr] \\
&&\qquad\leq e^{\rho T}
\int_{T-1}^{T} \,du\, \mathbb{P} \bigl[ e^{- \int_{- \infty}^{\infty}
(L_u^y)^2 \,dy} \mathbh{1}_{X_u \in[0,M]} \bigr]
\\
&&\qquad\leq e^{\rho} \int_{T-1}^{T} \,du\, \mathbb{P} \bigl[
e^{ \int_{ - \infty}^{\infty} [ - (L_u^y)^2 + \rho L_u^y ]
\,dy} \mathbh{1}_{X_u \in[0,M]} \bigr]
\end{eqnarray*}
by using the equality
\[
\int_{\mathbb{R}} \rho L_u^y \,dy = \rho u.
\]
By dominated convergence, Proposition \ref{prop1} is proven if we show that
%
%
\begin{equation} \label{xxxsss}
\int_{0}^{\infty} du\, \mathbb{P} \bigl[ e^{ \int_{ - \infty}^{\infty}
[ - (L_u^y)^2 + \rho L_u^y ] \,dy}
\mathbh{1}_{X_u \in[0,M]} \bigr] < \infty.
\end{equation}
In order to estimate the left-hand side of (\ref{xxxsss}), we need the
following lemma.
\begin{lemma} \label{leuridan}
For every positive and measurable function $G$ on $\mathbb{R} \times
\mathcal{C}(\mathbb{R}, \mathbb{R}_+)$
\[
\int_0^{\infty} \mathbb{P} [ G(X_u, L_u^{\bullet}) ] \,du =\int
_{\mathbb{R}} da \int_0^{\infty} dl\, \mathbb{E} [
G(a,Y_{l,a}^{\bullet}) ],
\]
where the law of the process $(Y_{l,a}^y)_{y \in\mathbb{R}}$ is
defined in the following way:
\begin{itemize}
\item for $a \geq0$, $(Y_{l,a}^{-y})_{y \geq0}$ is a squared Bessel
process of dimension zero, starting at~$l$;
\item for $a \geq0$, $(Y_{l,a}^{y})_{y \geq0}$ is an independent
inhomogeneous Markov process, which has the same infinitesimal
generator as a two-dimensional squared Bessel process for $y \in[0,a]$
and the same infinitesimal generator as a zero-dimensional squared
Bessel process for $y \geq a$;
\item for $a \leq0$, $(Y_{l,a}^{y})_{y \in\mathbb{R}}$ has the same
law as $(Y_{l,-a}^{-y})_{y \in\mathbb{R}}$.
\end{itemize}
\end{lemma}
\begin{pf}
For $a \geq0$, let $B$ be a standard Brownian motion,
$B^{(a)}$ an independent Brownian motion starting at $a$, and let us
denote by $(\tau_l)_{l \geq0}$ the inverse local time of $B$ at level
0, and $T^{(a)}_0$ the first time when $B^{(a)}$ reaches zero.

By \cite{l98} and \cite{by88}, for every process $(F_u)_{u \geq0}$ on
the space $\mathcal{C} (\mathbb{R}_+, \mathbb{R})$, which is
progressively measurable with respect to the filtration $(\mathcal
{F}_u)_{u \geq0}$,
%
%
\begin{equation} \label{0007}
\int_0^{\infty} du \,\mathbb{E} [F_u (B)] = \int_0^{\infty} dl \int
_{-\infty}^{\infty} da\, \mathbb{E} \bigl[F_{\tau_l + T^{(a)}_0} \bigl(Z^{(l,a)}
\bigr) \bigr],
\end{equation}
where $Z^{(l,a)}$ is a\vspace*{-2pt} process such that $Z^{(l,a)}_r = B_r$ for $r
\leq\tau_l$ and $Z^{(l,a)}_{\tau_l + T^{(a)}_0 - s} = B^{(a)}_s$ for
$s \leq T^{(a)}_0$.

By applying (\ref{0007}) to the process defined by $F_u(X) = G(X_u,
L_u^{\bullet})$, and by using Ray--Knight theorems, one obtains Lemma
\ref{leuridan}.
\end{pf}

An immediate application of this lemma is the following equality:
%
%
\begin{eqnarray} \label{nh28}
&&\int_{0}^{\infty} du\, \mathbb{P} \bigl[ e^{ \int_{ - \infty}^{\infty}
[ - (L_u^y)^2 + \rho L_u^y ] \,dy} \mathbh{1}_{X_u \in
[0,M]} \bigr]
\nonumber\\[-8pt]\\[-8pt]
&&\qquad= \int_{0}^{\infty} dl \int_{0}^{M} da\, \mathbb{E}
\bigl[ e^{ \int_{ - \infty}^{\infty} [ - (Y_{l,a}^y)^2 + \rho
Y_{l,a}^y ] \,dy}
\bigr].\nonumber
\end{eqnarray}
In order to majorize this expression, let us prove another result,
which is also used in the proof of Proposition \ref{prop2}.
\begin{lemma} \label{jj}
For all $l > 0$, $\mu\in\mathbb{R}$ and for all measurable functions
$g$ from $\mathbb{R}_+$ to $\mathbb{R}_+$, the following equality holds:
%
%
\begin{equation} \label{qwe}
\mathbb{E} \biggl[ e^{\int_{0}^{\infty} [ - (Y_{l,0}^y)^2 + \mu
Y_{l,0}^y ] \,dy} g \biggl(
\int_0^{\infty} Y_{l,0}^y \,dy \biggr) \biggr] = \int_0^{\infty}
K_l^{(\mu)} (v) g(v) \,dv,
\end{equation}
where $K_l^{(\mu)} (v)$ is defined by (\ref{kim28}).

In particular,
\[
\mathbb{E} \bigl[ e^{\int_{0}^{\infty} [ - (Y_{l,0}^y)^2 + \rho
Y_{l,0}^y ] \,dy} \bigr] = \bar{K}^{(\rho)}_l,
\]
where
\[
\bar{K}^{(\rho)}_l = \int_0^{\infty} K_l^{(\rho)} (v) \,dv.
\]
Moreover, $\bar{K}^{(\rho)}_l$ is bounded by a universal constant and
\[
\int_0^{\infty} \bar{K}^{(\rho)}_l \,dl < \infty.
\]
\end{lemma}
\begin{pf}
The process $Y_{l,0}$ is a local martingale with
bracket given, for $y \geq0$, by
\[
\langle Y_{l,0}, Y_{l,0} \rangle_y = 4 \int_0^y Y_{l,0}^x \,dx.
\]
Therefore,
\[
Y_{l,0}^y = 2 B_{\int_0^y Y_{l,0}^x \,dx}^{(l/2)}
\]
where $B^{(l/2)}$ is a Brownian motion starting at $l/2$. Moreover,
since $Y_{l,0}$ stays at zero when it hits 0, the hitting time of zero for
$B^{(l/2)}$ is $S = \int_0^{\infty} Y_{l,0}^x \,dx$. Hence, the change
of variable $s = \int_0^{y} Y_{l,0}^x \,dx$ gives
\[
\int_0^{\infty} (\mu- Y_{l,0}^y) Y_{l,0}^y \,dy = \int_0^S \bigl(\mu-
2 B_s^{(l/2)}\bigr) \,ds.
\]
Therefore, one has the equalities
\begin{eqnarray*}
&&\mathbb{E} \biggl[ e^{\int_{0}^{\infty} [ - (Y_{l,0}^y)^2 + \mu
Y_{l,0}^y ] \,dy} g \biggl(
\int_0^{\infty} Y_{l,0}^y \,dy \biggr) \biggr] \\
&&\qquad= \mathbb{E} \bigl[
e^{\int_{0}^{S} ( \mu-2 B^{(l/2)}_s) \,ds} g(S) \bigr]
\\
&&\qquad= \int_0^{\infty} e^{\mu v} g(v) \mathbb{E} [ e^{- \int
_{0}^{v} 2 B^{(l/2)}_s \,ds} | S=v ] \mathbb{P} [S \in dv].
\end{eqnarray*}
Now, this formula implies (\ref{qwe}), since the density at $v$ of the
law of $S$ is equal to $\alpha_l(v)$ and the law of $(B^{(l/2)}_s)_{s
\leq v}$, conditionally on $S=v$, is equal to the law of $V^{(l/2,v)}$
(see, e.g., \cite{gs79}).

It only remains to prove the integrability of $\bar{K}^{(\rho)}_l$. One
easily checks that
\[
\bar{K}^{(\rho)}_l = \mathbb{E} \bigl[ e^{\int_0^S (\rho- 2
B^{(l/2)}_s) \,ds } \bigr].
\]
Hence, if one sets
\[
f(x) = \operatorname{Ai}\bigl(2^{-1/3} (2x - \rho)\bigr),
\]
for the Airy function $\operatorname{Ai}$ [which is, up to a
multiplicative constant, the unique bounded solution of the differential
equation $\operatorname{Ai}'' (x) = x \operatorname{Ai}(x)$], the
process $N$ defined by
\[
N_t = f\bigl(B_t^{(l/2)}\bigr) \exp\biggl( \int_0^t \bigl( \rho- 2 B^{(l/2)}_s\bigr) \,ds
\biggr)
\]
is a local martingale.

Moreover, since $\rho$ is smaller than $-2^{1/3}$ times the largest
zero of Airy function, the function $f$ is strictly positive on
$\mathbb{R}_+$ and $N$ is positive.
By stopping $N$ at time~$S$, one gets a true martingale since $0 \leq
N_{t \wedge S} \leq\Vert f\Vert_{\infty} e^{\rho t}$, and
by Doob's stopping theorem and Fatou's lemma, one has
\[
\bar{K}^{(\rho)}_l \leq\frac{f(l/2)}{f(0)}.
\]
Since Airy function decays faster than exponentially at infinity, the
boundedness and the integrability of $\bar{K}^{(\rho)}_l$ are proven.
\end{pf}

It is now easy to prove that Lemma \ref{jj} implies Proposition \ref
{prop1}: by using this lemma, the definition of $Y_{l,a}$ and Markov
property at
level $a$, one can see that the left-hand side of (\ref{nh28}) is
equal to
%
%
\begin{equation} \label{bidule}
\int_{0}^{\infty} dl \int_{0}^{M} da\, \bar{K}^{(\rho)}_l \mathbb
{E} \bigl[e^{ \int_{0}^{a} [ - (Y_{l}^y)^2 + \rho Y_{l}^y ]}
\bar{K}^{(\rho)}_{Y_l^a} \bigr].
\end{equation}
Now, $\bar{K}^{\rho}_{Y_l^a} $ is uniformly bounded and $-x^2 + \rho x
\leq\frac{\rho^2}{4}$ for all $x \in\mathbb{R}$; hence, the quantity
(\ref{bidule})
is bounded by a constant times
\[
e^{M \rho^2 /4} \int_0^{M} da \int_0^{\infty} dl\, \bar{K}^{(\rho
)}_l,
\]
which is finite.

Hence, one has (\ref{xxxsss}), and finally Proposition \ref{prop1}.
\begin{remark*} Proposition \ref{prop1} remains true if one replaces
$\rho$ by any real $\rho'$ which is strictly smaller than $-2^{1/3}$
times the largest zero of Airy function (e.g., one can take
$\rho
' = 2.9$).
\end{remark*}

\section{\texorpdfstring{Proof of Proposition \protect\ref{prop2} (first part)}{Proof of Proposition 1.4 (first part)}}\label{sec4}

The purpose of this first part is to prove the following proposition,
which, in particular, gives another expression for the left-hand side
of (\ref{a}).
\begin{proposition} \label{palm}
For $u, v, t, l> 0$, let us define the quantities
%
%
\begin{equation} \label{20}
J_l (u,v) = \mathbb{E} \bigl[ e^{- 2 \int_0^u R_w^{(l/2)} dw} \chi
_v\bigl(2 R_u^{(l/2)}\bigr) \bigr],
\end{equation}
where $\chi_v$ is given by (\ref{guet}) and $(R_w^{(l/2)})_{w \geq0}$
is a Bessel process of dimension 2, starting at $l/2$;
%
%
\begin{equation} \label{21}
J_l(t) := \int_0^t J_l(t-v,v) \,dv
\end{equation}
and for all $t \in\mathbb{R}$,
\[
J_l^{(\rho)} (t) = e^{\rho t} \mathbh{1}_{t > 0} J_l(t).
\]
Then, there exists a subset $E$ of $\mathbb{R}_+$, such that the
complement of $E$ is Lebesgue-negligible, and for all $T \in E$,
%
%
\begin{eqnarray}\label{2008}
&&
e^{\rho T} \mathbb{P} \bigl[ e^{ - \int_{- \infty}^{\infty} [L_T^y +
f(y)]^2 \,dy} \mathbh{1}_{X_T \geq M} \bigr]
\nonumber\\[-8pt]\\[-8pt]
&&\qquad= \int_0^{\infty} dl\, \mathbb{E} \biggl[ e^{ \int_{-\infty}^{M} (
- [ Y_{l}^y + f(y) ]^2 + \rho Y_l^y ) \,dy}
J^{(\rho)}_{Y_{l}^M} \biggl( T - \int_{-\infty}^{M} Y_{l}^y \,dy
\biggr) \biggr].\nonumber
\end{eqnarray}
Moreover, for all measurable functions $h$ from $(\mathbb{R}_+)^2$ to
$\mathbb{R}_+$, and for all $l > 0$,
%
%
\begin{eqnarray}\label{666}
&&
\int_0^{\infty} db\, \mathbb{E} \biggl[ e^{- \int_0^{\infty}
(Y_{l,b}^y)^2 \,dy} h \biggl( \int_0^b Y_{l,b}^y \,dy , \int
_b^{\infty} Y_{l,b}^y \,dy \biggr) \biggr]
\nonumber\\[-8pt]\\[-8pt]
&&\qquad
= \int_{(\mathbb{R}_+)^2} h(u,v) J_l (u,v) \,du \,dv\nonumber
\end{eqnarray}
and for all measurable functions $g$ from $\mathbb{R}_+$ to
$\mathbb{R}_+$,
%
%
\begin{equation} \label{667}
\int_0^{\infty} db \,\mathbb{E} \biggl[ e^{- \int_0^{\infty}
(Y_{l,b}^y)^2 \,dy} g \biggl( \int_0^{\infty} Y_{l,b}^y \,dy
\biggr) \biggr] = \int_{0}^{\infty} g(t) J_l (t) \,dt.
\end{equation}
\end{proposition}
\begin{remark*} In Proposition \ref{palm}, it is natural to expect
that $E$ is empty, even if we do not need it to prove our main
result.
\end{remark*}
\begin{pf*}{Proof of Proposition \ref{palm}}
Equation (\ref{2008}) is a
consequence of (\ref{666}) and (\ref{667}); therefore, we begin our
proof by these two
equalities. By monotone class theorem, it is sufficient to
prove (\ref{666}) for functions $h$ of the form: $h(x,y) = h_1(x)
h_2(y)$, where $h_1$ and $h_2$ are measurable functions from $\mathbb
{R}_+$ to $\mathbb{R}_+$.

By Lemma \ref{jj}, for all $l > 0$,
\[
\mathbb{E} \biggl[ e^{- \int_0^{\infty} (Y_{l,0}^y)^2 \,dy} h_2
\biggl( \int_0^{\infty} Y_{l,0}^y \,dy \biggr) \biggr] = \int_0^{\infty} K_l
(v) h_2(v) \,dv.
\]
Hence, by applying Markov property to the process $Y_{l,b}$ at level
$b$,
%
%
\begin{eqnarray} \label{127}
&&\int_0^{\infty} db\, \mathbb{E} \biggl[ e^{- \int_0^{\infty}
(Y_{l,b}^y)^2 \,dy} h_1 \biggl( \int_0^b Y_{l,b}^y \,dy \biggr)
h_2 \biggl( \int_b^{\infty} Y_{l,b}^y \,dy \biggr) \biggr]
\nonumber\\
&&\qquad= \int_0^{\infty} db\, \mathbb{E} \biggl[ e^{- \int_0^{b} (Y_{l}^y)^2
\,dy} h_1 \biggl( \int_0^b Y_{l}^y \,dy \biggr) \int_0^{\infty}
K_{Y_{l}^b} (v) h_2(v) \,dv \biggr]
\\
&&\qquad
= \int_0^{\infty} dv\, h_2 (v) \mathbb{E} \biggl[ \int_0^{\infty} db\,
e^{- \int_0^{b} (Y_{l}^y)^2 \,dy} h_1 \biggl( \int_0^b Y_{l}^y
\,dy \biggr) K_{Y_{l}^b} (v) \biggr].\nonumber
\end{eqnarray}
Now, the function $\widetilde{y}$ from $\mathbb{R}_+$ to $\mathbb{R}_+$,
given by
\[
\widetilde{y}(s) = \inf\biggl\{ y \in\mathbb{R}_+, \int_0^y Y_l^{y'} \,dy' = s
\biggr\}
\]
is well defined, continuous, strictly increasing and tending to
infinity at infinity.

Hence, one can consider the process $(\widetilde{Q}_s^{(l)})_{s \geq0}$
such that
\[
\widetilde{Q}_s^{(l)} = Y_l^{\widetilde{y} (s)}.
\]
One has
\[
d \widetilde{y}(s) = \frac{ds}{Y_l^{\widetilde{y} (s)}} = \frac{ds}{\widetilde
{Q}_s^{(l)}},
\]
and the s.d.e.
\[
d \widetilde{Q}_s^{(l)} = 2 \sqrt{ Y_l^{\widetilde{y} (s)}} \,d \widehat
{B}_{\widetilde
{y}(s)} + 2 \,d\widetilde{y}(s) = 2 \,dB_s + \frac{2 \,ds}{\widetilde
{Q}_s^{(l)}},
\]
where $\widehat{B}$ and $B$ are Brownian motions: the processes $\widetilde
{Q}^{(m)}$ and $2 R^{(m/2)}$ have the same law.

By a change of variable in (\ref{127}) [$b= \widetilde{y} (s), y =
\widetilde{y} (u)$],
\begin{eqnarray*}
&&\int_0^{\infty} db\, \mathbb{E} \biggl[ e^{- \int_0^{\infty}
(Y_{l,b}^y)^2 \,dy} h_1 \biggl( \int_0^b Y_{l,b}^y \,dy \biggr)
h_2 \biggl( \int_b^{\infty} Y_{l,b}^y \,dy \biggr) \biggr]
\\
&&\qquad= \int_0^{\infty} dv\, h_2 (v) \mathbb{E} \biggl[ \int_0^{\infty} d
\widetilde{y}(s) e^{- \int_0^{s} (Y_{l}^ {\widetilde{y}(u)})^2 \,d \widetilde
{y} (u)} \,h_1 (s) K_{Y_{l}^{\widetilde{y} (s)}} (v) \biggr]
\\
&&\qquad= \int_0^{\infty} dv\, h_2 (v) \mathbb{E} \biggl[ \int_0^{\infty}
ds\,
h_1(s) e^{- 2 \int_0^{s} R^{(l/2)}_u \,du} \frac{K_{2
R^{(l/2)}_s} (v)}{2 R^{(l/2)}_s} \biggr],
\end{eqnarray*}
which implies (\ref{666}).

The equality (\ref{667}) is easily obtained by applying (\ref{666}) to
the function $h \dvtx (u,v) \rightarrow g(u+v)$.

Now, it remains to deduce (\ref{2008}) from (\ref{666}) and (\ref{667}).

For all measurable and positive functions $g$, one has, by Lemma
\ref{leuridan},
\begin{eqnarray*}
&&\int_0^{\infty} dT\, g(T) \mathbb{P} \bigl[ e^{- \int_{- \infty
}^{\infty} [L_T^y + f(y)]^2 \,dy} \mathbh{1}_{X_T \geq M} \bigr]
\\
&&\qquad= \int_0^{\infty} dl \int_M^{\infty} da \,\mathbb{E} \biggl[ e^{-
\int_{-\infty}^{\infty} [ Y_{l,a}^y + f(y) ]^2 \,dy} g
\biggl( \int_{-\infty}^{\infty} Y_{l,a}^y \,dy \biggr) \biggr].
\end{eqnarray*}
On the other hand, by applying Markov property (for the process
$Y_{l,a}$, at level $M$), and by using the fact that $f(y) = 0$ for $y
\geq M$, one obtains, for all positive and measurable functions $h_1$
and $h_2$,
%
%
\begin{eqnarray} \label{truc}\qquad
&&\int_0^{\infty} dl \int_M^{\infty} da\, \mathbb{E} \biggl[ e^{- \int
_{-\infty}^{\infty} [ Y_{l,a}^y + f(y) ]^2 \,dy} h_1
\biggl( \int_{-\infty}^{M} Y_{l,a}^y \,dy \biggr) h_2 \biggl( \int
_{M}^{\infty} Y_{l,a}^y \,dy \biggr) \biggr]
\nonumber\\
&&\qquad= \int_0^{\infty} dl \int_0^{\infty} db\, \mathbb{E} \biggl[ e^{-
\int_{-\infty}^{M} [ Y_{l}^y + f(y) ]^2 \,dy}
h_1 \biggl( \int_{-\infty}^{M} Y_{l}^y \,dy \biggr)\cdots
\\
&&\qquad\quad\hspace*{75.3pt}{}\times
\mathbb{E} \biggl[ e^{-
\int_{0}^{\infty} ( \widehat{Y}_{ Y_{l}^M, b}^y )^2 \,dy}
h_2 \biggl( \int_{0}^{\infty} \widehat{Y}_{ Y_{l}^M, b}^y \,dy
\biggr) \Big| Y_l^M \biggr] \biggr],\nonumber
\end{eqnarray}
where $\widehat{Y}_{ Y_{l}^M, b}$ is a process which has, conditionally
on $Y_{l}^M = l'$, the same law as $Y_{l',b}$.

Now, by putting the integral with respect to $db$ just before the
second expectation in the right-hand side of (\ref{truc}), and by
applying (\ref{667}) to $g=h_2$, one obtains that the left-hand side
of (\ref{truc}) is equal to
\[
\int_0^{\infty} dl\, \mathbb{E} \biggl[ e^{- \int_{-\infty}^{M}
[ Y_{l}^y + f(y) ]^2 \,dy} h_1 \biggl( \int_{-\infty}^{M}
Y_{l}^y \,dy \biggr) \int_0^{\infty} h_2(t) J_{Y_{l}^M} (t)
\,dt \biggr].
\]
Hence, by monotone class theorem, for all measurable functions $h$ from
$\mathbb{R}_+^2$ to $\mathbb{R}_+$,
\begin{eqnarray*}
&&\int_0^{\infty} dl \int_M^{\infty} da\, \mathbb{E} \biggl[ e^{- \int
_{-\infty}^{\infty} [ Y_{l,a}^y + f(y) ]^2 \,dy} h
\biggl( \int_{-\infty}^{M} Y_{l,a}^y \,dy , \int_{M}^{\infty} Y_{l,a}^y
\,dy \biggr) \biggr]
\\
&&\qquad
= \int_0^{\infty} dl\, \mathbb{E} \biggl[ e^{- \int_{-\infty}^{M}
[ Y_{l}^y + f(y) ]^2 \,dy}
\int_0^{\infty} h \biggl( \int_{-\infty}^{M} Y_{l}^y \,dy , t \biggr)
J_{Y_{l}^M} (t) \,dt \biggr].
\end{eqnarray*}
By applying this equality to the function $h \dvtx (u,v) \rightarrow
g(u+v)$, we obtain
\begin{eqnarray*}
&&
\int_0^{\infty} dT g(T) \mathbb{P} \bigl[ e^{- \int_{- \infty
}^{\infty} [L_T^y + f(y)]^2 \,dy} \mathbh{1}_{X_T \geq M} \bigr]
\\
&&\qquad= \int_0^{\infty} dT\, g(T) \int_0^{\infty} dl\, \mathbb{E} \biggl[
e^{- \int_{-\infty}^{M} [ Y_{l}^y + f(y) ]^2 \,dy}\cdots
\\
&&\qquad\quad\hspace*{101.7pt}{}\times J_{Y_{l}^M} \biggl(T - \int_{-\infty}^{M} Y_{l}^y \,dy
\biggr) \mathbh{1}_{\int_{-\infty}^{M} Y_{l}^y \,dy < T} \biggr].
\end{eqnarray*}
Since this equality is true for all $g$, there exists a subset $E$ of
$\mathbb{R}_+$, such that the complement of $E$ is Lebesgue-negligible,
and for all $T \in E$,
\begin{eqnarray*}
&&\mathbb{P} \bigl[ e^{- \int_{- \infty}^{\infty} [L_T^y + f(y)]^2 \,dy}
\mathbh{1}_{X_T \geq M} \bigr]
\\
&&\qquad= \int_0^{\infty} dl\, \mathbb{E} \biggl[ e^{- \int_{-\infty}^{M}
[ Y_{l}^y + f(y) ]^2 \,dy}
J_{Y_{l}^M} \biggl( T - \int_{-\infty}^{M} Y_{l}^y \,dy \biggr) \mathbh
{1}_{\int_{-\infty}^{M} Y_{l}^y \,dy < T} \biggr],
\end{eqnarray*}
which implies (\ref{2008}).
\end{pf*}

\section{Some estimates}\label{sec5}

In this section, we prove the following propositions, which give
estimates for the different quantities introduced earlier in this
paper.
In the sequel of this paper, $C$ denotes a universal and strictly
positive constant, which may change from line to line.
\begin{proposition} \label{81}
For all $l, v > 0$, $\mu\in\mathbb{R}$, one has the majorization
%
%
\begin{equation} \label{811}
K_l^{(\mu)} (v) \leq C l v^{-3/2} e^{(\mu-2.9) v - l^2/8v},
\end{equation}
where $K_l^{(\mu)} (v)$ is defined by (\ref{kim28}).
\end{proposition}
\begin{proposition} \label{82}
For all $M > 0$, the random variable $\int_0^M Y_0^y \,dy$ admits a
density $D_M$ with respect to Lebesgue measure, such that for all $u
\geq0$,
\[
\bigl(D_M \ast K_l^{(\rho)}\bigr) (u) \leq C_M e^{- \nu_M l},
\]
where $C_M, \nu_M > 0$ depend only on $M$.
\end{proposition}
\begin{proposition} \label{83}
For all $l, u, v > 0$
%
%
\begin{equation}\label{831}
J_l (u,v) \leq\frac{C e^{-2.9v}}{(u+v) \sqrt{v}} e^{-l^2/{8(u+v)}}
\end{equation}
and
%
%
\begin{equation} \label{832}
J_l (u,v) \leq\frac{C}{\sqrt{v}} e^{-2.8 v - \rho u},
\end{equation}
if $u \geq2$ [recall that $J_l(u,v)$ is defined by (\ref{20})].

Moreover, the function $\chi_v$ from $\mathbb{R}_+^*$ to $\mathbb{R}$
[recall that $\chi_v(l) = \frac{K_l(v)}{l}$], is in $L^{2}(\nu)$, and
for fixed $l, v > 0$ and $u$ going to infinity,
%
%
\begin{equation} \label{833}
e^{\rho u} J_l (u,v) \longrightarrow \langle \chi_v | e_0 \rangle
e_0 (l).
\end{equation}
\end{proposition}
\begin{proposition} \label{84}
For all $t > 0$
\[
e^{\rho t} J_l (t) \leq C \biggl( 1 + \frac{1}{\sqrt{t}} \biggr),
\]
where $J_l(t)$ is defined by (\ref{21}).

Moreover, for $l$ fixed and $t$ going to infinity
%
%
\begin{equation} \label{999999}
e^{\rho t} J_l (t) \longrightarrow K e_0(l),
\end{equation}
where $K$ is the universal constant defined in Proposition \ref{prop2}.
\end{proposition}
\begin{pf*}{Proof of Proposition \ref{81}}
For all $l \geq0$, the process $V^{(l/2,v)}$ is, by coupling,
stochastically larger than $V^{(l/2,v)}$. Therefore, by scaling
property,
\[
\mathbb{E} [ e^{-2 \int_0^v V_u^{(l/2,v)} \,du } ]
\leq\mathbb{E} [ e^{-2 \int_0^v V_u^{(0,v)} \,du } ] =
\mathbb{E} [ e^{-2 v^{3/2} \int_0^1 V_u^{(0,1)} \,du } ].
\]
Now, the Laplace transform of $\int_0^1 V_u^{(0,1)} \,du$ (the area
under a normalized Brownian excursion) is known (see, e.g.,
\cite{lou84});
one has, for $\lambda> 0$,
\[
\mathbb{E} [ e^{-\lambda\int_0^1 V_u^{(0,1)} \,du } ] =
\sqrt{2 \pi} \lambda\sum_{n=1}^{\infty} e^{-u_n (\lambda
^2/2)^{1/3}},
\]
where $-u_1 > -u_2 > -u_3 > \cdots$ are the (negative) zeros of the Airy
function.

Therefore,
\[
\mathbb{E} [ e^{-2 \int_0^v V_u^{(l/2,v)} \,du } ]
\leq\sqrt{8 \pi v^3 } \sum_{n=1}^{\infty} e^{-2^{1/3} u_n v}
\]
and
%
%
\begin{equation} \label{32}
\mathbb{E} [ e^{-2 \int_0^v V_u^{(l/2,v)} \,du } ]
\leq (C e^{-2.9v}) \sum_{n=1}^{\infty} e^{-(2^{1/3} u_n - 2.91)
v},
\end{equation}
since $v^{3/2}$ is dominated by $e^{0.01 v}$.

Now, $2^{1/3} u_1 > 2.91$; hence, for $v > 1$, the infinite sum in
(\ref{32}) is smaller than
\[
\sum_{n=1}^{\infty} e^{-(2^{1/3} u_n - 2.91)},
\]
which is finite, since $u_n$ grows sufficiently fast with $n$ (as $n^{2/3}$).
Consequently, for $v \geq1$,
\[
\mathbb{E} [ e^{-2 \int_0^v V_u^{(l/2,v)} \,du } ]
\leq C e^{-2.9v}.
\]
This majorization, which remains obviously true for $v \leq1$ if we choose
$C > e^{2.9}$, implies easily (\ref{811}).
\end{pf*}
\begin{pf*}{Proof of Proposition \ref{82}}
In \cite{bpy01}, the density
of the law of $\int_0^1 Y_0^y \,dy$ is explicitly given:
\[
D_1(x) = \pi\sum_{n=0}^{\infty} (-1)^n \biggl( n + \frac{1}{2}
\biggr) e^{- (n+{1}/{2} )^2 \pi^2 x /2}.
\]
This formula proves that $D_1$ is continuous on $\mathbb{R}_+^*$ and
that for $x \geq1$
%
%
\begin{equation} \label{key}
D_1 (x) \leq\pi e^{-\pi^2 (x-1)/8} \sum_{n=0}^{\infty} \biggl( n
+ \frac{1}{2} \biggr) e^{- (n+{1}/{2} )^2 \pi^2 /2} \leq
C e^{-x}.
\end{equation}
Moreover, $D_1$ satisfies the functional equation
\[
D_1 (x) = \biggl( \frac{2}{\pi x} \biggr)^{3/2} D_1 \biggl( \frac{4}{\pi
^2 x} \biggr),
\]
which proves that, for $x \leq\frac{4}{\pi^2}$,
\[
D_1 (x) \leq C x^{-3/2} e^{-4/\pi^2x} \leq C.
\]
This inequality and the continuity of $D_1$ imply that (\ref{key})
applies for all $x \in\mathbb{R}_+^*$.

By scaling property of squared Bessel processes, the density $D_M$
exists and one has
\[
D_M (x) =\frac{1}{M^2} D_1 \biggl( \frac{x}{M^2} \biggr),
\]
which implies
\[
D_M(x) \leq\frac{C}{M^2} e^{ - {x}/{M^2}}.
\]
Therefore, for all $u \geq0$,
\begin{eqnarray*}
\bigl(D_M \ast K_l^{(\rho)}\bigr) (u) & = & \int_0^u D_M(u-v) K_l^{(\rho)} (v)
\,dv \\
& \leq &\frac{C}{M^2} \int_0^u e^{ - ({u-v})/{M^2}}
l v^{-3/2} e^{(\rho-2.9) v - l^2/8v} \,dv \\
& \leq &\frac{C}{M^2}
e^{- ( 0.7 \wedge{1}/{M^2} ) u} \int_0^u l v^{-3/2}
e^{- l^2/8v} \,dv.
\end{eqnarray*}
The last inequality comes from the fact that $\rho-2.9 \leq0.7$,
which implies, for $0 \leq u \leq v$,
\[
- \frac{u-v}{M^2} + (\rho-2.9)v \leq- \biggl( \frac{u-v}{M^2} + 0.7 v
\biggr) \leq- \biggl( 0.7 \wedge\frac{1}{M^2} \biggr) u.
\]
Now, the integral $\int_0^u l v^{-3/2} e^{- l^2/8v} \,dv$ is
proportional to the probability that a Brownian motion starting at
$l/2$ reaches zero before time $u$.

Hence,
\[
\int_0^u l v^{-3/2} e^{- l^2/8v} \,dv \leq C e^{-l^2/8u},
\]
and finally
\begin{eqnarray*}
\bigl(D_M \ast K_l^{(\rho)}\bigr) (u) & \leq &\frac{C}{M^2}
e^{- 0.7u/(1+M^2) - l^2/8u} \\ & \leq &\frac{C}{M^2}
e^{- 2 \sqrt{( {0.7u}/({1+M^2}))(l^2/8u)}} \\ & \leq &
\frac{C}{M^2} e^{-{l}/({2(1+M)})},
\end{eqnarray*}
which proves Proposition \ref{82}.
\end{pf*}
\begin{pf*}{Proof of Proposition \ref{83}} By definition of $J_l(u,v)$,
one has
\[
J_l(u,v) \leq\mathbb{E} \bigl[ \chi_v \bigl(2 R_u^{(l/2)} \bigr)
\bigr].
\]
Now, the majorization (\ref{811}) implies
\[
\chi_v \bigl(2 R_u^{(l/2)} \bigr) \leq C v^{-3/2} e^{-2.9 v}
e^{- ( R_u^{(l/2)} )^2 / 2v}.
\]
By using the explicit expression of the Laplace transform of the
squared bidimensional Bessel process (see, e.g., \cite{ry91}),
one obtains
\[
J_l(u,v) \leq C v^{-3/2} e^{-2.9 v} \frac{v}{u+v}
e^{-l^2/8(u+v)},
\]
which implies (\ref{831}).

In order to prove (\ref{832}), let us consider, on the set of measurable
functions from $\mathbb{R}_+^*$ to $\mathbb{R}_+$, the semigroup of operators
$(\Phi^s)_{s \geq0}$ associated to the process $2R$
(twice a Bessel process of dimension 2), and the semigroup
$(\widetilde{\Phi}^s)_{s \geq0}$ associated to the same process,
killed at rate $x$ at level $x$.

For all positive and measurable functions $\psi$, and for all $l >
0$, one has
%
%
\begin{eqnarray}
\label{poix0}
[\Phi^s (\psi)](l) &=& \mathbb{E} \bigl[ \psi\bigl(2 R_s^{(l/2)}\bigr)
\bigr],
\\
\label{poix}
[\widetilde{\Phi}^s (\psi)](l) &=& \mathbb{E} \bigl[ e^{- 2 \int_0^s
R_u^{(l/2)} \,du} \psi\bigl(2 R_s^{(l/2)}\bigr) \bigr].
\end{eqnarray}
Now, let us observe that the measure $\nu$ on $\mathbb{R}_+^*$ is
reversible, and hence invariant by the semigroup of $2R$. Since, for
every measurable and positive function $\psi$,
\[
(\widetilde{\Phi}^s (\psi) )^2 \leq( \Phi^s (\psi)
)^2 \leq\Phi^s (\psi^2),
\]
one gets
%
%
\begin{equation} \label{poi}
\Vert\widetilde{\Phi}^s (\psi)\Vert_{L^2(\nu)}^2 \leq\int_{\mathbb{R}_+^*}
\psi^2 \,d \nu= \Vert\psi\Vert_{L^2(\nu)}^2.
\end{equation}
Inequality\vspace*{1pt} (\ref{poi}) proves that the semigroup
$(\widetilde{\Phi}^s)_{s \geq0}$ can be considered as a semigroup of
continuous linear operators on $L^{2} (\nu)$.

Moreover, the infinitesimal generator of $2R$, killed at rate
$x$ at level $x$, is the operator $\mathcal{K}$ defined at the
beginning of our paper. Hence, if $(e_n)_{n \in\mathbb{N}}$ is an
orthonormal basis of $L^2 (\nu)$ such that $e_n$ is an eigenvector of
$\mathcal{K}$,
corresponding to the eigenvalue $- \rho_n$ ($-\rho= - \rho_0 > -
\rho_1 \geq- \rho_2 \geq- \rho_3 \geq\cdots$), one has, for all $n$,
\[
\widetilde{\Phi}^s (e_n) = e^{- \rho_n s} e_n.
\]
Now, for all $\psi\in L^{2} (\nu)$, one has the representation
\[
\psi= \sum_{n \geq0} \langle \psi|e_n \rangle e_n,
\]
and, by linearity and continuity of $\widetilde{\Phi}^s$,
%
%
\begin{equation}\label{label}
\widetilde{\Phi}^s (\psi) = \sum_{n \geq0} e^{- \rho_n s}
\langle \psi|e_n \rangle e_n.
\end{equation}
In particular,
\[
\Vert \widetilde{\Phi}^s (\psi) \Vert _{L^{2}(\nu)}^2 =
\sum_{n \geq0} e^{-2 \rho_n s} (\langle \psi|e_n \rangle)^2 \leq e^{-2
\rho s}
\sum_{n \geq0} (\langle \psi|e_n \rangle)^2,
\]
which implies
%
%
\begin{equation}\label{38}
\Vert\widetilde{\Phi}^s (\psi)\Vert_{L^{2}(\nu)} \leq e^{- \rho s}
\Vert\psi\Vert_{L^2(\nu)}.
\end{equation}
Moreover, one has the equality
\[
e^{\rho s} \widetilde{\Phi}^s (\psi) - \langle \psi| e_0 \rangle
e_0 =
\sum_{n=1}^{\infty} e^{(\rho-\rho_n)s} \langle \psi| e_n \rangle
e_n,
\]
which implies that
\[
e^{\rho s} \widetilde{\Phi}^s (\psi)
\mathop{\longrightarrow}\limits_{s \rightarrow\infty}
\langle \psi| e_0 \rangle e_0
\]
in $L^2 (\nu)$.

Now, by definition
%
%
\begin{equation} \label{pira}
J_{l} (u,v) = (\widetilde{\Phi}^u (\chi_v) ) (l),
\end{equation}
where $\chi_v \in L^2 (\nu)$, by the majorization (\ref{811}).

Hence, in $L^2 (\nu)$,
\[
e^{\rho u} J_{\bullet} (u,v)
\mathop{\longrightarrow}\limits_{s \rightarrow\infty}
\langle \chi_v | e_0 \rangle e_0,
\]
where $ J_{\bullet} (u,v)$ is the function defined by
\[
(J_{\bullet} (u,v) ) (l) = J_{l} (u,v).
\]
In order to prove the corresponding pointwise convergence [which is
(\ref{833})], let us observe that for all $\psi\in L^{2} (\nu)$, $l
> 0$,
%
%
\begin{eqnarray} \label{79}
|[ \widetilde{\Phi}^1(\psi)] (l)|
&\leq&
\mathbb{E} \bigl[\bigl|\psi\bigl(2
R_1^{(l/2)}\bigr)\bigr|\bigr] \leq\bigl(\mathbb{E} \bigl[\bigl(\psi\bigl(2 R_1^{(l/2)}\bigr)\bigr)^2\bigr]
\bigr)^{1/2} \nonumber\\
&\leq& \biggl[ \int_0^{\infty} p_1^{(2)}
\biggl( \frac{l}{2},x \biggr) (\psi(2x))^2 \,dx \biggr] ^{1/2}
\nonumber\\[-8pt]\\[-8pt]
&\leq&\biggl[ \int_0^{\infty} x (\psi(2x))^2 \,dx
\biggr] ^{1/2} \nonumber\\
&\leq&\Vert\psi\Vert_{L^{2} (\nu)}.\nonumber
\end{eqnarray}
Here, we use the majorization $p_1^{(2)} (x,y) \leq y$, which comes
from the fact
that the transition densities of a bidimensional Brownian motion are
uniformly bounded by $1/2 \pi$ at time 1.

By (\ref{79}), one has for $s > 1$, $l > 0$, $\psi\in L^2(\nu)$,
\begin{eqnarray*}
|e^{\rho s} ( \widetilde{\Phi}^s (\psi) ) (l) - \langle
\psi| e_0 \rangle e_0 (l)| & = & \bigl| \bigl(e^{\rho s} \widetilde
{\Phi}^s (\psi) - \langle \psi| e_0 \rangle e_0 \bigr) (l)\bigr|
\\
& = & e^{\rho} \bigl| \bigl(\widetilde{\Phi}^1 \bigl(e^{\rho(s-1)}
\widetilde{\Phi}^{s-1} (\psi) - \langle \psi| e_0 \rangle e_0
\bigr) \bigr) (l)\bigr| \\
& \leq & e^{\rho} \bigl\Vert e^{\rho(s-1)} \widetilde{\Phi}^{s-1} (\psi) -
\langle \psi| e_0 \rangle e_0 \bigr\Vert _{L^2(\nu)} \\
&\mathop{\longrightarrow}\limits_{s \rightarrow\infty} &
0.
\end{eqnarray*}
By applying this convergence to $\chi_v$, one obtains the pointwise
(and in fact uniform) convergence (\ref{833}).

Now, it remains to prove (\ref{832}).

For $s \geq1$, $l > 0$, by (\ref{38}) and (\ref{79}),
%
%
\begin{eqnarray} \label{uuu}
|[\widetilde{\Phi}^s (\psi)] (l) | &=& |[\widetilde{\Phi}^1
(\widetilde
{\Phi}^{s-1} (\psi))](l)| \leq\Vert\Phi^{s-1} (\psi)\Vert_{L^2 (\nu)}
\nonumber\\[-8pt]\\[-8pt]
&\leq&
e^{- \rho(s-1)} \Vert\psi\Vert_{L^2 (\nu)}.\nonumber
\end{eqnarray}
By (\ref{pira}), and by semigroup property of $\widetilde{\Phi}$, one
has for $u >1$, $l, v > 0$,
%
%
\begin{equation} \label{vvv}
J_l (u,v) = [\widetilde{\Phi}^{u-1} (J_{\bullet} (1,v)) ](l).
\end{equation}
Moreover, by (\ref{831}),
\[
J_{l} (1,v) \leq C e^{-2.9v} v^{-1/2} e^{-l^2/8(1+v)},
\]
which implies
\begin{eqnarray*}
\Vert J_{\bullet} (1,v)\Vert_{L^2 (\nu)}^2 & \leq & C e^{-5.8v} \int
_0^{\infty
} \frac{l}{v} e^{-l^2/4(1+v)} \,dl \\ & \leq & C e^{-5.8v}
\frac{1+v}{v}
\\ & \leq &\frac{C}{v} e^{-5.6v}
\end{eqnarray*}
and
\[
\Vert J_{\bullet} (1,v)\Vert _{L^2 (\nu)} \leq\frac{C}{\sqrt{v}} e^{-2.8v}.
\]
For $u \geq2$, we can combine (\ref{uuu}) and (\ref{vvv}), and we obtain
\[
J_l (u,v) \leq C e^{- \rho(u-2)} \Vert J_{\bullet} (1,v)\Vert_{L^2 (\nu)}
\leq\frac{C}{\sqrt{v}} e^{- \rho u -2.8v},
\]
which is (\ref{832}).

The proof of Proposition \ref{83} is now complete.
\end{pf*}
\begin{pf*}{Proof of Proposition \ref{84}} Let us split the integral
corresponding to $J_l(t)$ into two parts:
\begin{eqnarray*}
A(t) &:=& \int_0^{(t-2)_+} J_l (t-v,v) \,dv,
\\
B(t) &:=& \int_{(t-2)_+}^t J_l (t-v,v) \,dv.
\end{eqnarray*}
One has
\[
e^{\rho t} A(t) = \int_0^{\infty} \mathbh{1}_{t-v \geq2} J_l (t-v,v)
e^{\rho t} \,dv;
\]
where, by (\ref{833}),
\[
\mathbh{1}_{t-v \geq2} J_l (t-v,v) e^{\rho t} \longrightarrow
e^{\rho v} \langle \chi_v | e_0 \rangle e_0(l),
\]
for $l, v$ fixed and $t$ tending to infinity.
Moreover, by (\ref{832}) and the fact that $\rho-2.8 < -0.6$,
\[
\mathbh{1}_{t-v \geq2} J_l (t-v,v) e^{\rho t} \leq\frac{C}{\sqrt
{v}} e^{-2.8 v - \rho(t-v) + \rho t} \leq
\frac{C}{\sqrt{v}} e^{-0.6 v},
\]
which is integrable on $\mathbb{R}_+$.

Hence,
%
%
\begin{equation} \label{x1}
e^{\rho t} A(t) \leq C,
\end{equation}
and by dominated convergence,
%
%
\begin{equation} \label{x2}
e^{\rho t} A(t)
\mathop{\longrightarrow}\limits_{t \rightarrow\infty}
K e_0 (l),
\end{equation}
where $l$ is fixed and
\[
K = \int_0^{\infty} e^{\rho v} \langle \chi_v |
e_0 \rangle \,dv
\]
is the constant defined in Proposition \ref{prop2}. The majorization
(\ref{x1}) implies that $K$ is
necessarily finite.

The integral $B(t)$ can be estimated in the following way: by (\ref{831}),
%
%
\begin{eqnarray}\label{x3}
e^{\rho t} B(t) & \leq & e^{\rho t} \int_{(t-2)_+}^{t} \frac{C
e^{-2.9v}}{t \sqrt{v}} e^{-l^2/8t} \,dv \nonumber\\ & \leq &
\frac{C e^{-0.7 t}}{t} e^{-l^2/8t} \int_0^t
\frac{dv}{\sqrt{v}} \nonumber\\[-8pt]\\[-8pt] & \leq &\frac{C e^{-0.7 t}}
{\sqrt{t}} e^{-l^2/8t} \nonumber\\ & \leq &\frac{C}
{\sqrt{t}}.\nonumber
\end{eqnarray}
Proposition \ref{84} is a consequence of (\ref{x1}), (\ref{x2}) and
(\ref{x3}).
\end{pf*}

The estimates given in this section are used in the second part of
the proof of Proposition \ref{prop2}.

\section{\texorpdfstring{Proof of Proposition \protect\ref{prop2} (second part)}{Proof of
Proposition 1.4 (second part)}}\label{sec6}

In this section, we estimate the right-hand side of (\ref{2008}), which
is, by Proposition \ref{palm}, equal to the left-hand side of (\ref
{a}) in
Proposition \ref{prop2}.

We need the following lemma.
\begin{lemma} \label{2007}
For all $M > 0$, and all functions $f, g$ from $\mathbb{R}_+$ to
$\mathbb{R}_+$,
\[
\int_0^{\infty} dl\, \mathbb{E} \biggl[ e^{ \int_{-\infty}^{M}
(- [ Y_{l}^y + f(y) ]^2 + \rho Y_l^y ) \,dy} g
\biggl( \int_{-\infty}^{M} Y_{l}^y \,dy \biggr) \biggr] \leq C'_M \int
_0^{\infty} g,
\]
where $C'_M>0$ is finite and depends only on $M$.
\end{lemma}
\begin{pf}
One has the following majorization:
%
%
\begin{eqnarray} \label{caipiri}
&&\int_0^{\infty} dl\, \mathbb{E} \biggl[ e^{ \int_{-\infty}^{M}
(- [ Y_{l}^y + f(y) ]^2 + \rho Y_l^y )
\,dy} g \biggl( \int_{-\infty}^{M} Y_{l}^y \,dy \biggr) \biggr]
\nonumber\\[-8pt]\\[-8pt]
&&\qquad\leq e^{M \rho^2 /4} \int_0^{\infty} dl\, \mathbb{E} \biggl[ e^{ \int
_{-\infty}^{0}
[ - ( Y_{l}^y )^2 + \rho Y_l^y ]
\,dy} g \biggl( \int_{-\infty}^{M} Y_{l}^y \,dy \biggr) \biggr]\nonumber
\end{eqnarray}
since $f$ is nonnegative and $- x^2 + \rho x \leq\rho^2/4$ for all
$x \in\mathbb{R}$.

Now, for all positive and measurable functions $h_1$ and $h_2$,
\begin{eqnarray*}
&&\int_0^{\infty} dl\, \mathbb{E} \biggl[ e^{ \int_{-\infty}^{0}
[- ( Y_{l}^y )^2 + \rho Y_l^y ] \,dy} h_1 \biggl(
\int_{-\infty}^{0} Y_{l}^y \,dy \biggr) h_2 \biggl( \int_{0}^{M}
Y_{l}^y \,dy \biggr) \biggr]
\\
&&\qquad= \int_0^{\infty} dl\, \mathbb{E} \biggl[ e^{ \int_{0}^{\infty}
[- ( Y_{l,0}^y )^2 + \rho Y_{l,0}^y ] \,dy} h_1
\biggl( \int_{0}^{\infty} Y_{l,0}^y \,dy \biggr) \biggr]
\mathbb{E} \biggl[ h_2 \biggl( \int_{0}^{M} Y_{l}^y \,dy \biggr) \biggr]
\\
&&\qquad= \int_0^{\infty} dl \int_0^{\infty} K_l^{(\rho)} (v) h_1(v)
\,dv\, \mathbb{E} \biggl[
h_2 \biggl( \int_{0}^{M} Y_{l}^y \,dy \biggr) \biggr],
\end{eqnarray*}
by Lemma \ref{jj}.

By additivity properties of squared Bessel processes, the law of $\int
_{0}^{M} Y_{l}^y \,dy $ is the convolution of the law $\sigma_1$ of
$\int_0^M Y_{l,0}^y \,dy$ and the law $\sigma_2$ of $\int_0^M Y_{0}^y
\,dy$.

Since by Proposition \ref{82}, $\sigma_2$ has the density $D_M$ with
respect to Lebesgue measure, we have the equality
\[
\mathbb{E} \biggl[
h_2 \biggl( \int_{0}^{M} Y_{l}^y \,dy \biggr) \biggr] = \int_0^{\infty}
dt\, h_2 (t) \int_0^t \sigma_1 (du)
D_M (t-u),
\]
which implies
\begin{eqnarray*}
&&\int_0^{\infty} dl\, \mathbb{E} \biggl[ e^{ \int_{-\infty}^{0}
[- ( Y_{l}^y )^2 + \rho Y_l^y ] \,dy} h_1 \biggl(
\int_{-\infty}^{0} Y_{l}^y \,dy \biggr) h_2 \biggl( \int_{0}^{M}
Y_{l}^y \,dy \biggr) \biggr]
\\
&&\qquad= \int_0^{\infty} dl \int_0^{\infty} K_l^{(\rho)} (v) h_1(v) \,dv
\int_0^{\infty} dt\, h_2 (t) \int_0^t \sigma_1 (du)
D_M (t-u).
\end{eqnarray*}
By monotone class theorem and easy computations, for all positive and
measurable functions $g$,
\begin{eqnarray*}
&&\int_0^{\infty} dl\, \mathbb{E} \biggl[ e^{ \int_{-\infty}^{0}
[- ( Y_{l}^y )^2 + \rho Y_l^y ] \,dy} g \biggl( \int
_{-\infty}^{M} Y_{l}^y \,dy \biggr) \biggr]
\\
&&\qquad= \int_0^{\infty} dt\, g(t) \int_0^{\infty} dl \int_0^t \sigma_1
(du) \bigl( K_l^{(\rho)} \ast D_M\bigr) (t-u).
\end{eqnarray*}
Now, by Proposition \ref{82},
\[
\bigl(K_l^{(\rho)} \ast D_M\bigr) (t-u) \leq C_M e^{- \nu_M l}
\]
and
\[
\int_0^t \sigma_1 (du) \leq1,
\]
since $\sigma_1$ is a probability measure.

Hence,
%
%
\begin{equation} \label{vi}
\int_0^{\infty} dl\, \mathbb{E} \biggl[ e^{ \int_{-\infty}^{0}
[- ( Y_{l}^y )^2 + \rho Y_l^y ] \,dy} g \biggl( \int
_{-\infty}^{M} Y_{l}^y \,dy \biggr) \biggr]
\leq\frac{C_M}{\nu_M} \int_0^{\infty} g.
\end{equation}
The majorizations (\ref{caipiri}) and (\ref{vi}) imply Lemma \ref
{2007}.
\end{pf}

After proving this lemma, let us take $T \in E$ and $\varepsilon> 0$; by
splitting the right-hand side of (\ref{2008}) into two parts, we obtain
\[
e^{\rho T} \mathbb{P} \bigl[ e^{ - \int_{- \infty}^{\infty} [L_T^y +
f(y)]^2 \,dy} \mathbh{1}_{X_T \geq M} \bigr] = I_{1,\varepsilon} +
I_{2, \varepsilon},
\]
where
\begin{eqnarray*}
I_{1, \varepsilon} &=& \int_0^{\infty} dl\, \mathbb{E} \biggl[ e^{ \int
_{-\infty}^{M} ( - [ Y_{l}^y + f(y) ]^2 + \rho Y_l^y
) \,dy}\cdots
\\
&&\hspace*{41.41pt}{}\times J^{(\rho)}_{Y_{l}^M} \biggl( T - \int_{-\infty}^{M}
Y_{l}^y \,dy \biggr) \mathbh{1}_{\int_{-\infty}^{M} Y_{l}^y \,dy
\notin[T- \varepsilon,T] } \biggr]
\end{eqnarray*}
and
\begin{eqnarray*}
I_{2, \varepsilon} &=& \int_0^{\infty} dl\, \mathbb{E} \biggl[ e^{ \int
_{-\infty}^{M} ( - [ Y_{l}^y + f(y) ]^2 + \rho Y_l^y
) \,dy}\cdots
\\
&&\hspace*{41.17pt}{}\times J^{(\rho)}_{Y_{l}^M} \biggl( T - \int_{-\infty}^{M}
Y_{l}^y \,dy \biggr) \mathbh{1}_{\int_{-\infty}^{M} Y_{l}^y \,dy \in
[T- \varepsilon,T] } \biggr].
\end{eqnarray*}
By Proposition \ref{84},
\[
J^{(\rho)}_{Y_{l}^M} \biggl( T - \int_{-\infty}^{M} Y_{l}^y \,dy
\biggr) \mathbh{1}_{\int_{-\infty}^{M} Y_{l}^y \,dy \notin[T- \varepsilon
,T] }
\mathop{\longrightarrow}\limits_{T \rightarrow\infty} K e_0 (Y_{l}^M)
\]
and
\[
J^{(\rho)}_{Y_{l}^M} \biggl( T - \int_{-\infty}^{M} Y_{l}^y \,dy
\biggr) \mathbh{1}_{\int_{-\infty}^{M} Y_{l}^y \,dy \notin[T- \varepsilon
,T] }
\leq C \biggl( 1 + \frac{1}{\sqrt{\varepsilon}} \biggr).
\]
Since
\begin{eqnarray*}
&&\int_0^{\infty} dl\, \mathbb{E} \bigl[ e^{ \int_{-\infty}^{M} (
- [ Y_{l}^y + f(y) ]^2 + \rho Y_l^y ) \,dy} \bigr]
\\
&&\qquad\leq e^{M \rho^2/4} \int_0^{\infty} dl\, \mathbb
{E} \bigl[ e^{ \int_{0}^{\infty} [ - ( Y_{l,0}^y )^2 +
\rho Y_{l,0}^y ] \,dy} \bigr]
\\
&&\qquad= e^{M \rho^2/4} \int_0^{\infty} \bar{K}^{(\rho)}_l \,dl < \infty,
\end{eqnarray*}
one obtains
%
%
\begin{eqnarray} \label{1283}
I_{1, \varepsilon} &\mathop{\longrightarrow}\limits_{T \rightarrow
\infty}& K
\int_0^{\infty} dl\, \mathbb{E} \bigl[ e^{ \int_{-\infty}^{M} ( -
[ Y_{l}^y + f(y) ]^2 + \rho Y_l^y ) \,dy} e_0
(Y_{l}^M) \bigr] \nonumber\\[-8pt]\\[-8pt]
&=&
K A_{+}^{1,M} (f) < \infty,\nonumber
\end{eqnarray}
by dominated convergence.

On the other hand, by Proposition \ref{84},
\begin{eqnarray*}
I_{2, \varepsilon} &\leq& C \int_0^{\infty} dl\, \mathbb{E} \biggl[ e^{
\int_{-\infty}^{M} ( - [ Y_{l}^y + f(y) ]^2 + \rho
Y_l^y ) \,dy}\cdots
\\
&&\hspace*{51.7pt}{}\times\biggl(1 + \biggl(T - \int_{-\infty}^{M} Y_{l}^y \,dy
\biggr)^{-1/2} \biggr) \mathbh{1}_{\int_{-\infty}^{M} Y_{l}^y \,dy \in
[T- \varepsilon,T] } \biggr],
\end{eqnarray*}
and by applying Lemma \ref{2007} to the function $g \dvtx t \rightarrow(1+
(T-t)^{-1/2}) \mathbh{1}_{t \in[T- \varepsilon, T]}$,
%
%
\begin{equation}\label{1284}
I_{2, \varepsilon} \leq C C'_M \bigl(\varepsilon+ \sqrt{\varepsilon}\bigr).
\end{equation}
Therefore, by combining (\ref{1283}) and (\ref{1284}),
\[
\mathop{\lim\sup}_{ T \in E, T \rightarrow\infty} \bigl|e^{\rho T}
\mathbb{P} \bigl[ e^{ - \int_{- \infty}^{\infty} [L_T^y + f(y)]^2
\,dy} \mathbh{1}_{X_T \geq M} \bigr] - K A_{+}^{1,M} (f) \bigr| \leq C
C'_M \bigl(\varepsilon+ \sqrt{\varepsilon}\bigr),
\]
and by taking $\varepsilon\rightarrow0$,
%
%
\begin{equation} \label{t4}
e^{\rho T} \mathbb{P} \bigl[ e^{ - \int_{- \infty}^{\infty} [L_T^y +
f(y)]^2 \,dy} \mathbh{1}_{X_T \geq M} \bigr]
\mathop{\longrightarrow}\limits_{T \in E, T \rightarrow\infty}
K A_{+}^{1,M} (f).
\end{equation}
Now, let us prove the continuity, with respect to $T$, of the left-hand
side of (\ref{t4}).

If $T_0 \in\mathbb{R}_+$ and $T \leq T_0 + 1$ tends to $T_0$, then
$\mathbb{P}$-almost surely, $L_T^y$ tends to $L_{T_0}^y$ and $L_{T}^y
\leq L_{T_0 + 1}^y$ for all $y \in\mathbb{R}$.

Since $y \rightarrow L_{T_0 + 1}^y + f(y)$ is square-integrable, by
dominated convergence,
\[
\int_{- \infty}^{\infty} [L_T^y + f(y)]^2 \,dy
\mathop{\longrightarrow}\limits_{T \rightarrow T_0}
\int_{- \infty}^{\infty} [L_{T_0}^y + f(y)]^2
\,dy.
\]
Another application of dominated convergence gives
%
%
\begin{equation} \label{1010}
\bigl| \mathbb{P} \bigl[ e^{ - \int_{- \infty}^{\infty} [L_T^y + f(y)]^2
\,dy} \mathbh{1}_{X_{T_0} \geq M} \bigr]
- \mathbb{P} \bigl[ e^{ - \int_{- \infty}^{\infty} [L_{T_0}^y +
f(y)]^2 \,dy} \mathbh{1}_{X_{T_0} \geq M} \bigr]\bigr|
\mathop{\longrightarrow}\limits_{T \rightarrow T_0}
0.\hspace*{-25pt}
\end{equation}
Moreover,
\[
\mathbb{P} [ |\mathbh{1}_{X_{T_0} \geq M} - \mathbh{1}_{X_{T} \geq
M} | ] \leq
\mathbb{P} \bigl[ \exists t \in[T_0, T], X_t = M \bigr]
\mathop{\longrightarrow}\limits_{T \rightarrow T_0}
\mathbb{P} [
X_{T_0} = M ] = 0,
\]
which implies
%
%
\begin{equation} \label{1011}
\bigl| \mathbb{P} \bigl[ e^{ - \int_{- \infty}^{\infty} [L_T^y + f(y)]^2
\,dy} \mathbh{1}_{X_{T} \geq M} \bigr]
- \mathbb{P} \bigl[ e^{ - \int_{- \infty}^{\infty} [L_{T}^y + f(y)]^2
\,dy} \mathbh{1}_{X_{T_0} \geq M} \bigr]\bigr|
\mathop{\longrightarrow}\limits_{T \rightarrow T_0}
0.\hspace*{-25pt}
\end{equation}
The convergences (\ref{1010}) and (\ref{1011}) imply the continuity of
\[
T \longrightarrow e^{\rho T} \mathbb{P} \bigl[ e^{ - \int_{- \infty
}^{\infty} [L_T^y + f(y)]^2 \,dy} \mathbh{1}_{X_T \geq M} \bigr].
\]
Since $E$ is dense in $\mathbb{R}_+$, we can remove the condition $T
\in E$ in (\ref{t4}), which completes the proof of Proposition
\ref{prop2}.

\section{A conjecture about the behavior of $\mathbb{Q}^{\beta}$}\label{sec7}

In this paper, we have proven that one can construct a probability
measure corresponding to the one-dimensional Edwards' model, for
polymers of infinite length.

Moreover, there is an explicit expression for this probability
$\mathbb{Q}^{\beta}$.

Now, the most natural question one can ask is the following: what is
the behavior of the canonical process $X$ under $\mathbb{Q}^{\beta}$?

At this moment, we are not able to answer this question, which seems to
be very difficult, because of the complicated form of the density
$D_s^{\beta}$ of $\mathbb{Q}^{\beta}_{|\mathcal{F}_s}$,
with respect to $\mathbb{P}_{| \mathcal{F}_s}$.

However, it seems to be reasonable to expect that $X_T$ has a
ballistic behavior, as in the case of Edwards' model on
$[0,T]$; one can also expect a central-limit
theorem.

Therefore, we can state the following conjecture.
\begin{conjecture*} Under $\mathbb{Q}^{\beta}$, the process $X$ is
transient, and
\[
\mathbb{Q}^{\beta} \bigl( X_t
\mathop{\longrightarrow}\limits_{t \rightarrow\infty}
+ \infty\bigr) = \mathbb{Q}^{\beta} \bigl(
X_t
\mathop{\longrightarrow}\limits_{t \rightarrow\infty}
- \infty\bigr) = 1/2.
\]
Moreover, there exist universal positive constants $a$ and $\sigma$
such that
\[
\frac{|X_t|}{t}
\mathop{\longrightarrow}\limits_{t \rightarrow\infty}
a \beta^{1/3}
\]
a.s., and such that the random variable
\[
\frac{|X_t|- a \beta^{1/3} t }{\sqrt{t}}
\]
converges in law to a centered Gaussian variable of variance
$\sigma^2$ (the factor $\beta^{1/3}$ comes from the Brownian
scaling).
\end{conjecture*}

It is possible that the constants in these convergences are the same
as in \cite{hhk97}, despite the fact that we don't have any argument to
support this. It can also be
interesting to study some large deviation results for the canonical
process under~$\mathbb{Q}^{\beta}$,
and to compare them with the results given in \cite{hhk03}. On the
other hand, if the proof of the conjecture above
is too hard to obtain, it is perhaps less
difficult\vspace*{1pt} to prove, by using Ray--Knight theorems, some properties of
the total local times
$(L_{\infty}^y)_{y \in\mathbb{R}}$ of $X$, which are expected to be
finite because of the transience of $X$.

%

%
\printaddresses

\end{document}